\documentclass[12pt]{article}

\usepackage{amsmath}
\usepackage{amsfonts}
\usepackage{amssymb}
\usepackage{graphicx}

\newtheorem{proposition}{Proposition}[section]
\newtheorem{theorem}{Theorem}[section]
\newtheorem{lemma}[theorem]{Lemma}
\newtheorem{corollary}[theorem]{Corollary}
\newtheorem{remark}[theorem]{Remark}
\newtheorem{example}[theorem]{Example}

\newtheorem{definition}{Definition}

\def\phi{{\varphi}}

\DeclareSymbolFont{AMSb}{U}{msb}{m}{n}
\DeclareMathSymbol{\N}{\mathbin}{AMSb}{"4E}
\DeclareMathSymbol{\Z}{\mathbin}{AMSb}{"5A}
\DeclareMathSymbol{\R}{\mathbin}{AMSb}{"52}
\DeclareMathSymbol{\Q}{\mathbin}{AMSb}{"51}
\DeclareMathSymbol{\I}{\mathbin}{AMSb}{"49}
\DeclareMathSymbol{\C}{\mathbin}{AMSb}{"43}

\begin{document}

\title{Conductivity imaging from one interior measurement in the presence of perfectly conducting and insulating inclusions }

\author{{Amir Moradifam\footnote{Department of Mathematics, University of Toronto, Toronto, Ontario, Canada M5S 2E4. E-mail: amir@math.toronto.edu. The author is supported by a MITACS Postdoctoral Fellowship. }
\qquad Adrian Nachman\footnote{Department of Mathematics and the
Edward S. Rogers Sr. Department of Electrical and Computer
Engineering, University of Toronto, Toronto, Ontario, Canada M5S
2E4. E-mail: nachman@math.toronto.edu.}\qquad Alexandru
Tamasan\footnote{Department of Mathematics, University of Central
Florida, Orlando, FL, USA. E-mail: tamasan@math.ucf.edu. The work of this author was supported  by the NSF grant DMS-0905799.}}}
\date{\today}

\smallbreak \maketitle

\begin{abstract}
We consider the problem of recovering an isotropic conductivity
outside some perfectly conducting or insulating inclusions from the
interior measurement of the magnitude of one current density field
$|J|$. We prove that the conductivity outside the inclusions, and
the shape and position of the perfectly conducting and insulating
inclusions are uniquely determined (except in an exceptional case) by the magnitude of the current
generated by imposing a given boundary voltage.  We have found an extension of the notion of admissibility to the 
case of possible presence of perfectly conducting and insulating inclusions. This also makes it possible to extend the results on uniqueness 
of the minimizers of the least gradient problem
$F(u)=\int_{\Omega}a|\nabla u|$ with $u|_{\partial \Omega}=f$ to cases where $u$ has flat regions (is constant on open sets).

\end{abstract}
\maketitle

\section{Introduction}
This paper considers the inverse problem of determining an isotropic
electrical conductivity $\sigma$ from one measurement of the
magnitude of the current density field $|J|$ generated inside the
domain $\Omega$ while imposing the voltage $f$ at the boundary.
Extending the existing work, the problem here allows for some
perfectly conducting and insulating inclusions be embedded in
$\Omega$ away from the boundary. The
 domain $\Omega\subset \mathbb{R}^n$, $n\geq2$, is assumed bounded,
open and with a connected Lipschitz boundary.

The problem considered in this paper is modelled by two physical principles: the Maxwell
model of the electromagnetic field at very low frequency, and a
magnetic resonance technique to image current densities pioneered in \cite{joy89} and
\cite{joy}. Employment of dual physical models is a fairly new trend in
quantitative imaging which seeks better accuracy and resolution of
the reconstructed images, compared to the methods based on
just one physical principle. For recent progress in such hybrid
imaging methods in conductivity imaging we refer to
\cite{arridge06},
\cite{ABCTF},\cite{scherzer},\cite{ACKK},\cite{balSchotland},\cite{balUhlmann},
\cite{wang}, \cite{kuchmentKunyanski}, and the review articles \cite{B} and \cite{NTT11}.

Inspired by \cite{joy89} and \cite{joy}, two subclasses of conductivity imaging methods have been developed: the ones which use interior
knowledge of the current density field, and the ones that use the
measurement of only one component of the magnetic field, known as
Magnetic Resonance Electric Impedance Tomography (see
\cite{B_z2003},\cite{B_z2004},\cite{B_z2005},\cite{B_z2006},\cite{B_z2008},\cite{B_z2007},
\cite{B_z2010} for work in this direction). The problem considered
here belongs to the former subclass. The idea of using the current
density field to image electrical conductivity appeared first in
\cite{zhang}. In \cite{ider98} a perturbation method recovered the
conductivity in the linearized case. Using the fact that $J$ is
normal to equipotential lines, the method in \cite{jeun-rock}
recovered two dimensional isotropic conductivities. In \cite{seo}
the problem is reduced to the Neumann problem for the 1-Laplacian,
and the examples of non-uniqueness and non-existence for this
degenerate elliptic problem show that knowledge of the applied
current at the boundary together with the magnitude of current
density field inside is insufficient data to determine the
conductivity. Instead, the ``$J$- substitution'' algorithm based on
knowledge of the magnitude of two current density fields has been
proposed; see also \cite{seoIEEE02} and \cite{seo03}. The idea of
using two currents goes back to \cite{scott}; in \cite{lionheart}
the problem is reduced to a first order system of PDEs and several
numerical reconstructions based on solving this system are proposed.
In independent work in \cite{cdii04}, and respectively \cite{lee}, a
simple formula recovers $\nabla\ln(\sigma)$ at each point in a
region where two transversal current density vectors have been
measured; see also \cite{cdii08} for careful experimental validation of this formula.

In \cite{NTT07} a reconstruction method which uses the interior
knowledge of the magnitude of just one current density field $|J|$ has
been proposed. This method relies on the fact that, in the absence
of singularities, equipotential sets are minimal surfaces in the
metric $g=|J|^{2/(n-1)} I$ conformal to the Euclidean metric. In
\cite{NTT10} it is shown that the equipotential surfaces are
minimizers for the area functional
\begin{equation}\label{massfunctional}
\mathcal{A}(\Sigma)=\int_{\Sigma}|J|dS,
\end{equation}
where $dS$ is the induced Euclidean surface measure. (Note that
$\mathcal{A}(\Sigma)$ is the area of $\Sigma$ in the Riemannian
metric $g$ described above.) Moreover, in \cite{NTT08} it is shown that the voltage potential $u$ is a
minimizer of the functional
\begin{equation}\label{eu_lag}
\int_{\Omega}|J|\cdot|\nabla{v}|dx,
\end{equation}
subject to $v\in W^{1,1}(\Omega)$ with $v=f$ at the boundary
$\partial\Omega$, and that $u$ is the unique minimizer among $v\in
W^{1,1}(\Omega)$ with $|\nabla v|>0$ a.e. in $\Omega$ and $v=f$ at
the boundary. One can determine $u$, and hence $\sigma$ by a minimization algorithm. A
structural stability result for the minimization of the functional
in (\ref{eu_lag}) can be found in \cite{nashedTa10}. Formally, the
Euler-Lagrange equation for the non-smooth functional in
(\ref{eu_lag}) is the generalized 1-Laplacian. This is in contrast
with the work in \cite{ABCTF}, \cite{ACKK} and \cite{scherzer},
where the conductivity imaging from interior data leads to the generalized
0-Laplacian.

Partial reconstruction from incomplete data results
are available for planar domains \cite{NTT10}: If $|J|$ is known throughout $\Omega$, but $f$ is
only known on parts of the boundary. More precisely, if some
interval $(\alpha,\beta)$ of boundary voltages is twice contained in
the known values of $f$, then one can recover the conductivity in the
subregion
\begin{equation}\label{imaged}
\Omega_{\alpha,\beta}:=\{x\in\overline\Omega:\alpha<u(x)<\beta\}.
\end{equation}
In fact  $|J|$ need only be known in a subregion
$\tilde{\Omega}$ which contains
regions of the type (\ref{imaged}) for unknown values $\alpha$'s and
$\beta$'s. The method in \cite{NTT10} determines from the data if
$\tilde\Omega$ contains regions of the type (\ref{imaged}), and, if
so, recovers all the (maximal) intervals $(\alpha,\beta)$, their
corresponding $\Omega_{\alpha,\beta}$ and the conductivity therein.

In this paper we are interested in imaging an isotropic conductivity
$\sigma$  from the magnitude of one current density field in the presence of  perfectly conducting and insulating
inclusions. We shall
prove that the conductivity outside the inclusions, and the shape
and position of the perfectly conducting and insulating inclusions
are uniquely determined (except in an exceptional case, see Remark \ref{rem1}) by the magnitude of the current generated by
imposing a given boundary voltage. We also establish a connection
between the above problem and the uniqueness of the minimizers of
weighted least gradient problem $F(u)=\int_{\Omega}a|\nabla u|$ with
$u|_{\partial \Omega}=f$.

Unlike the results in \cite{NTT07}, \cite{NTT08}, and \cite{NTT10}
that have been proven under the assumption that the interior data
$|J|>0$ a.e. in $\Omega$, the results presented in this paper allow for $|J|\equiv 0$ in open subsets of $\Omega$. In the following section we present and
discuss our main results.

\section{Main results}

Let $U$ be an open subset of $\Omega$ with
$\overline{U}\subset\Omega$ to model the \emph{perfectly conducting
inclusions}, $V$ be an open subset of $\Omega$ with
$\overline{V}\subset\Omega$ to model the \emph{insulating
inclusions}, and let $\chi_U$ and $\chi_V$ be their corresponding
characteristic functions. Note that $U$ and $V$ may have more than one connected component. We assume
$\overline{U}\cap\overline{V}=\emptyset$,  $\Omega\setminus
\overline{U\cup V}$ is connected, and the boundaries $\partial U$,
$\partial V$ are piecewise $C^{1,\alpha}$. Let $\sigma_1\in
L^\infty(U)$, and  $\sigma\in L^\infty(\Omega\setminus
\overline{U\cup V})$ be bounded away from zero. For  $k>0$ consider
the conductivity problem
\begin{equation}\left\{ \begin{array}{ll}
\nabla\cdot((\chi_U(k\sigma_1-\sigma)+\sigma)\nabla u)=0, \ \ \hbox{in} \ \  \Omega \setminus \overline{V}\\
\frac{\partial u}{\partial \nu}=0 \ \  \hbox{on}\ \ \partial V, \\
 u|_{\partial\Omega}=f.
\end{array} \right.
\end{equation}

The perfectly conducting inclusions occur in the limiting case $k\to
\infty$. The limiting solution is the unique solution to the
problem:
\begin{equation}\left\{ \begin{array}{ll}
\nabla\cdot \sigma \nabla u_0=0,&\mbox{in}\,\Omega\setminus\overline{U\cup V},\\
\nabla u_0=0, &\mbox{in} \ \ U,\\
u_0|_+=u_0|_-,&\mbox{on}\ \ \partial (U\cup V),\\
\int_{\partial U_j}\sigma\frac{\partial u_0}{\partial \nu}|_{+}ds=0,&j=1,2,...,\\
\frac{\partial u_0}{\partial \nu}|_{+}=0,&\mbox{on}\;\partial V,\\
u_0|_{\partial\Omega}=f,
\end{array}\label{pde_inclusions} \right.
\end{equation}
 (see the Appendix for more details), where  $U=\cup_{j=1}^\infty U_j$ is a  partition of $U$ into connected
components.


For Lipschitz continuous conductivities in any dimension $n\geq 2$,
or for essentially bounded conductivities in two dimensions, the
solutions of the conductivity equation satisfy the unique
continuation property (see, \cite{bersJohnSchechter} and references
therein). Consequently the insulated (and possibly perfectly
conducting) inclusions are the only open sets on which the interior
data $|J|$ vanishes identically. However, in three dimensions or
higher it is possible to have a H\"{o}lder continuous $\sigma$ and
boundary data $f$ that yield $u\equiv constant$ in a proper open
subset $W\subsetneq \Omega$, see \cite{plis,martio}. We call such
regions $W$ {\em singular inclusions}. On the other hand Ohm's law
need not hold inside perfect conductors: the current $J$ inside
perfectly conducting inclusions $U$ is not necessarily zero while
$\nabla u\equiv 0$ in $U$ (\cite{Aetall}, \cite{lp}).

The measured data for our inverse problem is the non-negative function $a=|J(x)|$ in $\Omega$, the magnitude of the current density field $J$
induced by imposing a voltage $f$ at the boundary $\partial \Omega$.  We have $\nabla \cdot J=0$. In the perfectly  conducting inclusion $U$ we will not rely on the Ohm's law; we will use the condition (\ref{perfect-cond}) and the transmission  condition $J_{-}\cdot \nu=J_{+}\cdot \nu$ across the boundary of $\partial U$ (see the Appendix). Indeed we have found an extension of the notion admissibility of \cite{NTT08} which will be crucial in allowing us to treat the case of perfectly conducting and insulating inclusions considered here. In a different direction, this also makes it possible to extend results on uniqueness of minimizers of weighted least gradient problems as discussed later in this section.

To formulate our results, we first need to introduce a notion of \emph{admissibility}.

\begin{definition} \label{def} A pair of functions $(f,a)\in H^{1/2}(\partial \Omega)\times
L^{2}(\Omega)$ is called \emph{admissible} if the following
conditions hold:

(i) There exist two disjoint open sets $U,V \subset \Omega$ (possibly empty)
and a function $\sigma \in L^{\infty}(\Omega\setminus (U\cup V))$
bounded away from zero such that $\Omega \setminus (\overline{U\cup V})$ is connected and 
\begin{eqnarray*}\label{rel}
\left\{ \begin{array}{ll}
a=|\sigma \nabla u_\sigma| \ \ \hbox{in}\ \ \Omega \setminus (\overline{U\cup V}),\\
a=0 \ \ \hbox{in} \ \ V,
\end{array} \right.
\end{eqnarray*}
where $u_\sigma \in H^{1}(\Omega)$ is the weak solution of
(\ref{pde_inclusions}).

(ii) The following holds

\begin{equation}\label{perfect-cond}
\inf_{u\in W^{1,1}(U)} \left( \int_{U} a|\nabla u|-\int_{\partial U}\sigma \frac{\partial u_{\sigma}}{\partial \nu}|_{+} u \right )=0,
\end{equation}
where $\nu$ is the unit normal vector field on $\partial U$ pointing outside $U$.

(iii) The set of zeroes of  the function $a$ outside $\overline{U}$
can be partitioned as follows
\begin{equation}\label{zeroes }
\{x\in \Omega:\, a(x)=0\}\cap (\Omega\backslash \overline{U})= V\cup
\overline{W} \cup \Gamma,
\end{equation}
where $W$ is an open set (possibly empty) , $\Gamma$ is a
Lebesgue-negligible set, and $\overline{\Gamma}$ has empty interior.

We call $\sigma$ a \emph{generating
conductivity} and $u_\sigma$ the \emph{corresponding potential}.

\end{definition}

Since for $u=constant$,
\[\int_{U_j} a|\nabla u|-\int_{\partial U_j}\sigma \frac{\partial u_{\sigma}}{\partial \nu}|_{+} u=0, \]
we have
\[\inf_{u\in W^{1,1}(U_j)} \left( \int_{U_j} a|\nabla u|-\int_{\partial U_j}\sigma \frac{\partial u_{\sigma}}{\partial \nu}|_{+} u \right )\leq 0.\]
Hence the condition $(\ref{perfect-cond})$ holds if and only if 
\[\inf_{u\in W^{1,1}(U_j)} \left( \int_{U_j} a|\nabla u|-\int_{\partial U_j}\sigma \frac{\partial u_{\sigma}}{\partial \nu}|_{+} u \right )=0,\]
for all connected components $U_j$ of $U$.

We first note that any physical data $(f,a)$  naturally satisfies the first two conditions i) and ii) in the above definition.  Indeed if  $a=|J|$ where $\nabla \cdot J=0$ in $\Omega$, then for any $u \in W^{1,1}(U)$ we have
\begin{eqnarray*}
\int_{U} a|\nabla u|-\int_{\partial U}\sigma \frac{\partial u_{\sigma}}{\partial \nu} u &=&\int_{U} |J| |\nabla u|-\int_{\partial U}\sigma \frac{\partial u_{\sigma}}{\partial \nu} u \\
&\geq& \int_{U} J\cdot \nabla u-\int_{\partial U}\sigma \frac{\partial u_{\sigma}}{\partial \nu} u \\
&=&\int_{\partial U} J_{-} \cdot\nu u-\int_{\partial U}\sigma \frac{\partial u_{\sigma}}{\partial \nu} u\\
&=&\int_{\partial U} J_{-}\cdot \nu u-\int_{\partial U}J_{+}\cdot \nu u=0.
\end{eqnarray*}
Also by fourth equation in (\ref{pde_inclusions})
\[\int_{U} a|\nabla u|-\int_{\partial U}\sigma \frac{\partial u_{\sigma}}{\partial \nu} u=0,\]
for any constant function $u$ in $U$. Hence ii) holds for physical data $(f,a)$. The first condition i) also obviously holds for physical data $(f,a)$.  We have added condition (iii) for technical reasons. Even though it is not always satisfied, this condition is very general, at least for physical applications.

On the other hand if 
\[\int _{U}\sigma \frac{\partial u_\sigma}{\partial \nu} \neq 0,\]
then 
\[E(u)=\int_{U} a|\nabla u|-\int_{\partial U}\sigma \frac{\partial u_{\sigma}}{\partial \nu} u,\]
is not invariant under adding or subtracting constant and therefore 
\[ \inf _{u\in W^{1,1}(U)} \left( \int_{U} a|\nabla u|-\int_{\partial U}\sigma \frac{\partial u_{\sigma}}{\partial \nu} u \right)=\infty.\]
Thus we have the following proposition about condition (\ref{perfect-cond}). 

\begin{proposition}
Let $a \in L^{\infty} (\Omega)$ and $U$ be an open subset of $\Omega$. Then 
\begin{itemize}
\item If $a\geq |J|$ in $U$ for some $J$ with $\nabla \cdot J \equiv 0$ in $U$ and $J_{-}=\sigma \frac{\partial u_\sigma}{\partial \nu}|_{+}$ on $\partial U$, then the condition (\ref{perfect-cond}) in Definition \ref{def} holds. 
\item If the the condition (\ref{perfect-cond}) in Definition \ref{def} holds, then 
\[\int _{U}\sigma \frac{\partial u_\sigma}{\partial \nu} =0. \]
\end{itemize}
\end{proposition}
We can now state one of our main uniqueness results.

\begin{theorem} \label{unique}
Let $\Omega \subset \R^n$, $n\geq 2$, be a domain with connected
Lipschitz boundary and let $(f,|J|) \in C^{1,\alpha} (\partial
\Omega) \times L^2(\Omega)$ be an admissible pair generated by some
unknown  conductivity $\sigma \in C^{\alpha}(\Omega \backslash (\overline{U\cup
V}))$, where $U$ and $V$ are open sets as described in
Definition \ref{def}. Then the potential $u_\sigma$ is a minimizer of the problem
\begin{equation}\label{min_prob}
u=\hbox{argmin} \{\int_{\Omega}|J||\nabla v|: v \in W^{1,1}(\Omega),
\ \ v|_{\partial \Omega}=f\},
\end{equation}
and if $u$ is another minimizer of the above problem, then
$u=u_\sigma$ in
\[\Omega \backslash \{x\in \Omega: \ \ |J|=0\}.\]
Moreover the set of zeros of $|J|$ and $|\nabla u_\sigma |$ can be decomposed
as follows
\[\{x\in \Omega: \ \ |J|=0\}\cup \{x\in \Omega: \ \ \nabla u_\sigma=0\}=:Z\cup \Gamma,\]
where $Z$ is an open set and $\Gamma$ has measure zero and
\[Z=U\cup V\cup W.\]
Consequently $\sigma=\frac{|J|}{|\nabla u_\sigma|} \in
L^{\infty}(\Omega\setminus \overline {Z})$ is the unique
$C^{\alpha}(\Omega\setminus \overline {Z})$-conductivity outside $Z$
for which $|J|$ is the magnitude of the current density corresponding to the voltage $f$ at the boundary.
\end{theorem}

\begin{remark} \label{rem1}The above theorem allows us to identify the potential $u=u_\sigma$
and the conductivity $\sigma$ outside the open set $Z=U \cup V \cup
W$. There are number of ways to determine if an open connected
component $O$ of $Z$ is a perfectly conducting inclusion, an
insulating inclusion, or a singular inclusion:
\begin{itemize}
\item If $\nabla u\equiv 0$ in $O$ and $|J|(x)\neq 0$ for
some $x \in O$, then $O$ is a perfectly conducting inclusion.
\item If $|J|\equiv 0$ in $O$ and $ u\not\equiv constant$ on $\partial O$,
then $O$ is an insulating inclusion.
\item If $J\equiv0$ in $O$, $u=constant$ on $\partial O$, and $J$ is not
$C^{\alpha}$ at $x$ for some $x \in O$, then $O$ is either an
insulating inclusion or a perfectly conducting inclusion.
\item If $J\equiv0$, $u=constant$ on $\partial O$, and $J\in C^{\alpha}(\partial
O)$, then the knowledge of the magnitude of the current $|J|$ (and even the full vector field $J$) is not
enough to determine the type of the inclusion $O$.
\end{itemize}
\end{remark}

\begin{remark} On can compare the forward problem (\ref{pde_inclusions}) with the minimization problem (\ref{min_prob}) to see that second, third, fourth, and fifth condition in the forward problem (\ref{pde_inclusions}) do not appear in the problem (\ref{min_prob}).  This means that all of the information about the location and shape of the inclusions is encoded in $|J|$. 
\end{remark}

Now we introduce an interesting connection between Theorem
\ref{unique} and the uniqueness of minimizers of weighted least
gradient problems. Indeed, Theorem \ref{unique} can also be applied
independently to prove uniqueness of the minimizers of the weighted
least gradient problem
\begin{equation} \label{lg}
u_0=\hbox{argmin}\{\int_{\Omega}a|\nabla u|, \ \ u\in W^{1,1}(\Omega), \ \ \hbox{and}\ \ u|_{\partial \Omega}=f\},
\end{equation}
in situations where thje minimizer has flat regions (is constant on open sets).
\begin{example}\label{ex}
For instance consider the following example \cite{sternberg_ziemer}.  Let  $D=\{x\in R^2: \ \ x^2+y^2<1\}$ be the unit disk and $f(x,y)=x^2-y^2$.  Consider the problem
\begin{equation}\label{sz}
u_0=\hbox{argmin}\{\int_{D}|\nabla u|, \ \ u\in W^{1,1}(D), \ \ \hbox{and}\ \ u|_{\partial D}=f\},
\end{equation}
which corresponds to $a\equiv|J|\equiv1$ in $D$.  We claim that $(1,x^2-y^2)$ is an admissible pair according to Definition \ref{def}.  To prove our claim we let
$U=(-\frac{1}{\sqrt{2}}, \frac{1}{\sqrt{2}} )\times (-\frac{1}{\sqrt{2}}, \frac{1}{\sqrt{2}})$ and $V=\emptyset$. Define

\begin{eqnarray*}
\sigma=\left\{ \begin{array}{ll}
\frac{1}{4|x|}, \ \ \text{if } \ \ |x|\geq \frac{1}{\sqrt{2}}, \ \ |y|\leq \frac{1}{\sqrt{2}},\\
\frac{1}{4|y|}, \ \ \text{if } \ \ |x|\leq \frac{1}{\sqrt{2}}, \ \ |y|\geq \frac{1}{\sqrt{2}},\\
\end{array} \right.
\end{eqnarray*}
and

\begin{eqnarray*}
u_\sigma=\left\{ \begin{array}{ll}
2x^2-1, \ \ \text{if } \ \ |x|\geq \frac{1}{\sqrt{2}}, \ \ |y|\leq \frac{1}{\sqrt{2}},\\
0, \hspace{0.5in}\ \hbox{if} \ \ (x,y)\in U,\\
1-2y^2, \ \ \text{if } \ \ |x|\leq \frac{1}{\sqrt{2}}, \ \ |y|\geq \frac{1}{\sqrt{2}}.\\
\end{array} \right.
\end{eqnarray*}

It is easy to see that $u_\sigma$ is the solution of (\ref{pde_inclusions}) and $|J|\equiv1\equiv \sigma|\nabla u_\sigma|$ on $\Omega \setminus \overline{U}$.  Hence (i) holds in the definition of admissibility, Definition \ref{def}.  The condition  (iii) also obviously holds.  It remains to show that (\ref{perfect-cond}) holds.
Define the vector field $J(x,y)$ in $U$ as follows
\begin{eqnarray*}
J(x,y)=\left\{ \begin{array}{ll}
-j, \ \ &\text{if } \ \ y\geq |x|,\\
j, \ \ &\text{if } \ \ -y\geq|x|,\\
i, \ \ &\text{if } \ \ x>|y|,\\
-i, \ \ &\text{if } \ \ -x>|y|,\\
\end{array} \right.
\end{eqnarray*}

Let 
\[U_0=\{(x,y)\in U | \ \ |x|\neq |y|\}=T_1\cup T_2 \cup T_3 \cup T_4,\]
where $T_i$, $1\leq i \leq 4$, are the four disjoint triangles in Figure 1.  Then $|J|=1$ in $U$, $J \in C^{\infty}(U_0)$ and we have
 \begin{eqnarray*}
\int_{U}|\nabla u|-\int_{\partial U}\sigma \frac{\partial u_\sigma}{\partial \nu}u &\geq& \int_{U_0}|J||\nabla u|-\int_{\partial U}\sigma \frac{\partial u_\sigma}{\partial \nu}u\\
&\geq& \int_{U_0}J\cdot \nabla u-\int_{\partial U}\sigma \frac{\partial u_\sigma}{\partial \nu}u\\
&=&\sum_{i=1}^{4} \int_{T_i}J\cdot \nabla u-\int_{\partial U}\sigma \frac{\partial u_\sigma}{\partial \nu}u\\
&=&\int_{\partial U} J\cdot \nu u-\int_{\partial U}\sigma \frac{\partial u_\sigma}{\partial \nu}u\\
&=&0,
\end{eqnarray*}
since $J\cdot \nu \equiv \sigma \frac{\partial u_\sigma}{\partial \nu}$
on $\partial U$. Thus the condition (\ref{perfect-cond}) holds
and $(1,x^2-y^2)$ is  admissible in the sense of Definition \ref{def}. It follows from Theorem
\ref{unique} that $u_\sigma$ is the unique minimizer of the problem
(\ref{sz}). 

\begin{figure}[h]
  \centering
  \includegraphics[height=60mm]{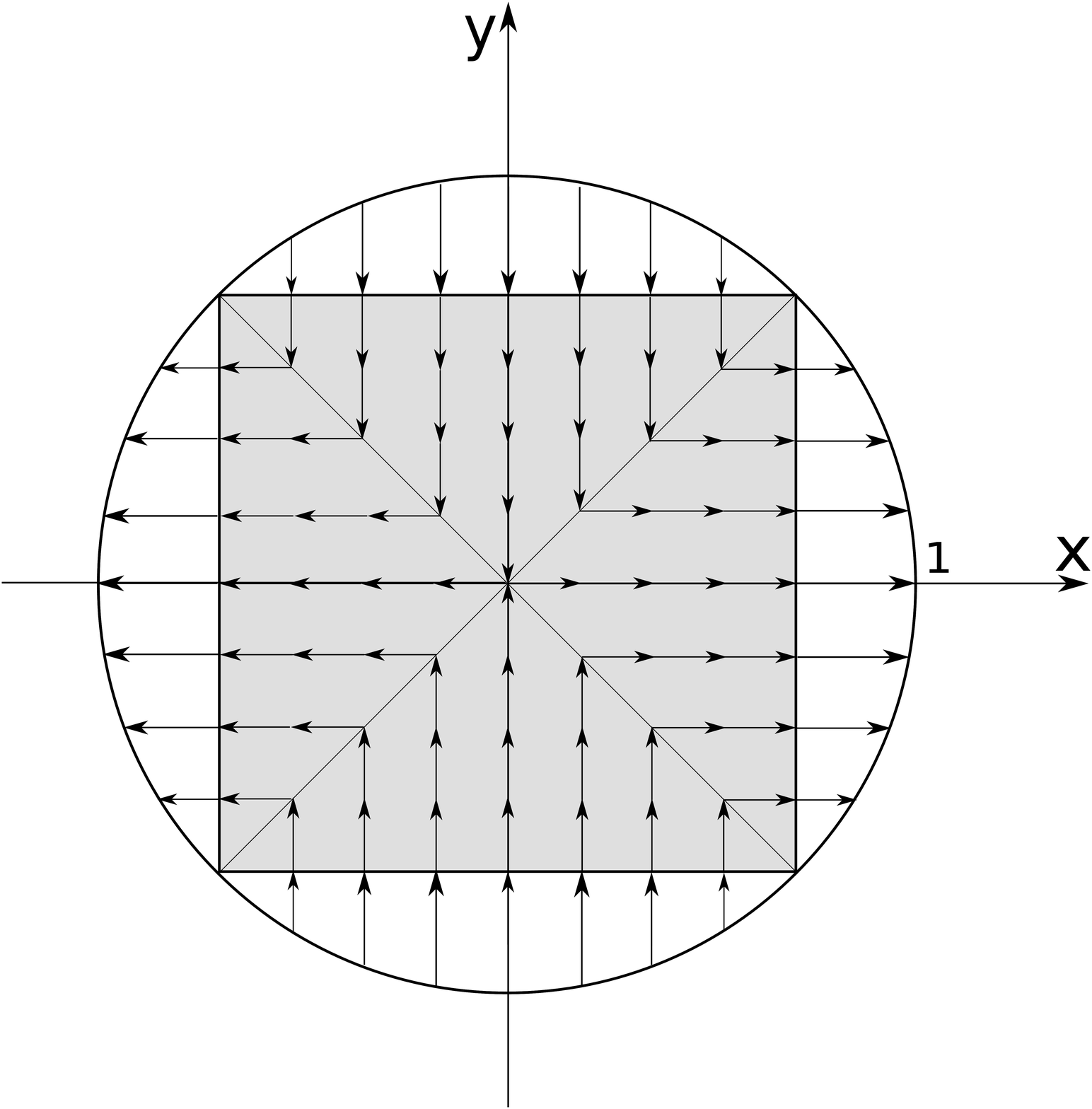}
  \caption{Current density vector field for Example \ref{ex} }
\end{figure}

\end{example}

The following theorem shows that the equipotential sets contained
entirely outside the conductive inclusions are area minimizers. We
describe a surface as the level set of a regular map $u$, while
competitors are described by level sets of some compact
perturbations of the regular map $u$.

\begin{theorem} \label{area.min}(Minimizing property of level sets). Let $\Omega \subset \R^n$, $n\geq 2$, be a domain with connected Lipschitz
boundary and let $(f,|J|) \in C^2
(\partial \Omega) \times L^2(\overline{\Omega})$ be an admissible
pair generated by some unknown $C^{1}$
conductivity. Then for every $v \in C^{2}(\overline{\Omega})$ with
$v=f$ on $\partial \Omega$ such that
\[\{x: |\nabla v|=0\}=Z_v\cap L_v, \ \ a(\overline{Z_v})=\{0\},\]
 where $Z_v$ is open and $L_v$ has Lebesgue
measure zero we have
\begin{equation}\label{ine}
{\cal A}(u^{-1}(\lambda))\leq {\cal A}(v^{-1}(\lambda)),
\end{equation}
for a.e. $\lambda \in \R$, where $\mathcal{A}$ is defined as
$(\ref{massfunctional})$.
\end{theorem}

The partial data result \cite[Theorem 3.4]{NTT10} also
recovers the conductivity in two dimensional subregions of type
(\ref{imaged}) assuming that $|J|> 0$ almost everywhere.
Below we show that, under the assumption the full vector field $J$
is known (not just its magnitude $|J|$), the partial reconstruction
result is valid in three  or higher dimensions. The result below can
be viewed as the extension of the results in \cite{jeun-rock} to
three or higher dimensional models.

\begin{theorem}\label{p.unique}
(Partial determination). Let $\Omega \subset \R^n$ ($n\geq 2$) be
simply connected. For $i=1,2$, let
$\sigma^i\in C^{\alpha}(\Omega \backslash \overline {U^i\cup V^i})$
be bounded away from zero, and $u_i$ satisfy
$(\ref{pde_inclusions})$, where $U^i$ and $V^i$ are open sets of
$\Omega$, and let
\begin{eqnarray*}
J_i=\left\{ \begin{array}{ll} \sigma^i \nabla u_i &\hbox{in}\ \
\Omega \backslash (U^i\cup V^i)\\
0 & \hbox{in} \ \  V^i
\end{array} \right.
\end{eqnarray*}
For $\alpha<\beta$ let
\begin{equation} \Omega_{\alpha,\beta}:=\{x \in \overline{\Omega}: \
\ \alpha <u_1<\beta \} \ \ \hbox{and} \ \
\Gamma:=\Omega_{\alpha,\beta} \cap
\partial \Omega.
\end{equation}
Assume that
\[\{x
\in \Omega\setminus \overline{U_1} : |J_1(x)|=0\}=V^1\cup W^1\cup \Gamma^1 ,\]
 where $W^1$ is open and $\Gamma^1$ has Lebesgue measure zero. Then
\begin{enumerate}
\item if $u_1|_{\Gamma}=u_2|_{\Gamma}$ and $J_1=J_2$ in
$\Omega $. Then $U^1\cap \Omega_{\alpha,\beta}= U^2\cap \Omega_{\alpha,\beta}$, $(W^1\cup V^1)\cap \Omega_{\alpha,\beta}=(W^2\cup V^2)\cap \Omega_{\alpha,\beta}$
\[u_1=u_2 \ \ \hbox{in} \ \ \Omega_{\alpha,\beta}\backslash V^1 \ \  \hbox{and} \ \ \sigma^1=\sigma^2 \ \ \hbox{in} \ \ \Omega_{\alpha,\beta}\setminus
\overline{U^1\cup V^1\cup W^1}.\]
 \item if $u_1|_{\Gamma}=u_2|_{\Gamma}$ and $J_1=J_2$ in $\Omega_{\alpha,\beta} $. Then
\begin{equation}
\{x\in \overline{\Omega}: \ \ \alpha<u_2(x)<\beta
\}=\Omega_{\alpha,\beta},
\end{equation}
$U^1\cap \Omega_{\alpha,\beta}= U^2\cap \Omega_{\alpha,\beta}$, $(W^1\cup V^1)\cap \Omega_{\alpha,\beta}=(W^2\cup V^2)\cap \Omega_{\alpha,\beta}$ and
\[u_1=u_2 \ \ \hbox{in} \ \ \Omega_{\alpha,\beta}\backslash V^1 \ \ \hbox{and} \ \ \sigma^1=\sigma^2 \ \ \hbox{in} \ \ \Omega_{\alpha,\beta} \setminus
\overline{U^1\cup V^1\cup W^1}.\]
\end{enumerate}
\end{theorem}

Similar to Theorem \ref{unique} we may determine if an open connected component $O$ of $U^1\cup V^1\cup W^1=U^2\cup V^2\cup W^2$ is a perfectly conducting, insulating, or singular inclusion (see Remark \ref{rem1}).

\section{Unique determination of the conductivity}
In this section we prove Theorems \ref{unique} and \ref{p.unique}.
The arguments extend those in \cite{NTT08} and \cite{NTT10}  by replacing the new admissibility condition. We
start with the following proposition.

\begin{proposition}
Let $\Omega \subset \R^n$, $n\geq 2$ be a domain and $(f,|J|)\in
H^{1/2}(\partial \Omega)\times L^2(\Omega)$. Then
\begin{enumerate}
\item Assume $(f,|J|)$ is admissible, say generated by some conductivity
$\sigma \in L^{\infty}(\Omega\backslash (\overline{U \cup V}))$
where $U$ and $V$ is described in Definition \ref{def} and $u_{0}$
is the corresponding voltage potential. Then $u_{0}$ is a minimizer
for $F(u)$ in (\ref{eu_lag}) over
\begin{equation}
A:=\{u\in H^{1}(\Omega): \ \ u|_{\Omega}=f\}.
\end{equation}
Moreover, if $f \in C^{1,\alpha}(\partial \Omega)$ and if the
generating conductivity $\sigma \in C^{\alpha}(\Omega \backslash
\overline {U \cup V})$, then the corresponding potential $u_{0} \in
C^{1,\alpha}(\Omega \backslash \overline {U \cup V})$ is a minimizer
of $F(u)$ over $A$.
\item Assume that the set of zeros of $a=|J|$ can be decomposed as follows
\[\{x\in \Omega: \ \ a(x)=0\}=V\cup \Gamma_1,\]
where $V$ is an open set and $\Gamma_1$ has measure zero. Suppose $u_0$ is a minimizer for $F(u)$ in (\ref{eu_lag}) over $A$ and the set of zeroes of $|\nabla u_0|$ can be decomposed as follows
\[\{x\in \Omega \setminus V : |\nabla u_0 |=0\}=\overline{U}\cup \Gamma_2,\]
where $U$ is an open set and $\overline{U \cup V}\subset \Omega$, and $\Gamma_2$ has measure
zero. If $U \cap V=\emptyset$ and  $|J|/|\nabla u_0| \in L^{\infty}(\Omega \backslash
(\overline{U\cup Z}))$ , then $(f,|J|)$ is admissible.
\end{enumerate}
\end{proposition}
{\bf Proof:} Assume $(f,|J|)$ is admissible and generated by some
conductivity $\sigma \in L^{\infty}_{+}(\Omega \backslash (U\cup
V))$. For any $u \in A$ we have
\begin{eqnarray*}
F(u)&=&\int_{\Omega \backslash \overline{(U\cup Z)}}\sigma |\nabla
u_0||\nabla u|dx+\int_{U}|J||\nabla u|dx \\
&\geq& \int_{\Omega \backslash \overline{(U\cup Z)}} \sigma \nabla
u_{0}.\nabla u+\int_{U}|J||\nabla u|dx \\
&=&\int_{\partial\Omega}\sigma \frac{\partial u_0}{\partial \nu}u ds-\int_{\partial V}\sigma \frac{\partial u_0}{\partial \nu}u ds-\int_{\partial U}\sigma \frac{\partial u_0}{\partial \nu}u ds+\int_{U}|J||\nabla u|dx\\
&=&\int_{\partial\Omega}\sigma \frac{\partial u_0}{\partial \nu}u ds-\int_{\partial U}\sigma \frac{\partial u_0}{\partial \nu}u ds+\int_{U}|J||\nabla u|dx\\
&\geq&\int_{\partial \Omega}\sigma \frac{\partial u_0}{\partial
\nu}f ds=F(u_0).
\end{eqnarray*}
where we have used the admissibility condition (\ref{perfect-cond}) and $\nu$ is the outer normal to the boundary of $\Omega$, $U$, and $V$. Hence $u_0$ is a minimizer of $F(u)$.

To prove 2) we note that by Lebesgue dominated convergence theorem,
the functional $F$ is Gateaux-differentiable at $u \in
H^{1}(\Omega)$ with $\frac{|J|}{|\nabla u|}\in L^{\infty}
_{+}(\Omega \backslash \overline{(U\cup V)})$. Since
\[F(u_0)=\int_{\Omega}|J||\nabla u_0|=\int_{\Omega \setminus \overline{U\cup V}}|J||\nabla u_0|,\]
 at a minimizer $u_0$ we have
\[F'(u_0)(\varphi)=\int_{\Omega \backslash \overline{U\cup V}}\frac{|J|}{|\nabla u_0|}\nabla u_0.\nabla \varphi dx=0,\]
for all $\varphi \in W_0^{1,1}(\Omega \setminus \overline {U})$.  Now let $\sigma=\frac{|J|}{|\nabla u_0|}$, then $\nabla.(\sigma \nabla u_0)=0$ in $\Omega \setminus  \overline{V\cup U}$. On the other hand we have
\[\int_{\Omega \backslash \overline{U\cup V}}\frac{|J|}{|\nabla u_0|}\nabla u_0.\nabla \varphi dx=\int_{\partial (U\cup V)}\sigma \frac{\partial u_0}{\partial \nu}\varphi dx=
\int_{\partial  V}\sigma \frac{\partial u_0}{\partial \nu}\varphi dx=0,\]
for all $\varphi \in W_0^{1,1}(\Omega \setminus \overline {U})$. Therefore $\frac{\partial u_0}{\partial \nu}=0$ on $\partial V$. Now let $ O$ be a connected component of $U$. Then for all $\varphi \in W^{1,1}_0(\Omega\setminus \overline{U\setminus O})$ with $\varphi\equiv 1$ in $O$ we have
\[\int_{\Omega \backslash \overline{U\cup V}}\frac{|J|}{|\nabla u|}\nabla u_0.\nabla \varphi dx=\int_{\partial (U\cup V)}\sigma \frac{\partial u_0}{\partial \nu}\varphi dx=
\int_{\partial  O}\sigma \frac{\partial u_0}{\partial \nu} dx=0.\]
This implies that $u_0$ is a solution of
(\ref{pde_inclusions}) (see the appendix for more details).

Moreover for every $u \in W^{1,1}_0(\Omega)$ with $u|_{\partial \Omega}=f$
\begin{eqnarray*}
\int_{\Omega}|J||\nabla u_0|  dx &\leq &\int_{\Omega\setminus \overline{V}}|J||\nabla u| dx\\
& =&\int_{U}|J||\nabla u| dx +\int_{\Omega \setminus \overline{U\cup V}}|J||\nabla u| dx \\
&=& \int_{U}|J||\nabla u| dx +\int_{\Omega \setminus \overline{U\cup V}}\sigma |\nabla u_0||\nabla u| \\
&=& \int_{U}|J||\nabla u| dx +\int_{\Omega \setminus \overline{U\cup V}}\sigma \nabla u_0.\nabla u \\
&=& \int_{U}|J||\nabla u| dx -\int_{\partial U}\sigma \frac{\partial u_0}{\partial \nu}u dx+\int_{\partial \Omega}\sigma \frac{\partial u_0}{\partial \nu}f dx.
\end{eqnarray*}
Since
\[\int_{\Omega}|J||\nabla u_0|  dx =\int_{\partial \Omega}\sigma \frac{\partial u_0}{\partial \nu} f dx,\]
the admissibility condition (\ref{perfect-cond}) follows from the above inequality. Thus $(|J|,f)$ is an admissible pair. \hfill $\Box$\\

Now we are ready to prove Theorem \ref{unique}. \\

{\bf Proof of Theorem \ref{unique}:} Assume $u_0$ is a solution of
(\ref{pde_inclusions}) that corresponds to the admissible pair $(f,
|J|)$. It is a direct consequence of the admissibility assumption that
\[\{x\in \Omega: \ \ |J|=0\}\cup \{x\in \Omega: \ \ \nabla u_\sigma=0\}=:Z\cup \Gamma,\]
where $Z$ is an open set and $\Gamma$ has measure zero and
\[Z=U\cup V\cup W.\]

Als well, since $\partial (U
\cup V) $ is piecewise $C^{1,\alpha}$,
\[ u_0 \in C^{1,\alpha}(\Omega \backslash \overline {U\cup V}) \cap C(\Omega \backslash \overline{U\cup V} \cup \partial
\Omega)\cap C^{1,\alpha}(\Omega \backslash \overline{U\cup V} \cup
T)\] for every $C^{1,\alpha}$ component of $\partial (U\cup V)$.

By our assumptions $|J|>0$ a.e. in $\Omega \setminus \overline{U\cup
V\cup W}$. Hence, equality in (\ref{rel}) yields $|\nabla u_0|>0$
a.e. on $\Omega \setminus \overline{U\cup V\cup W}$. Since $U\cup W$
is a disjoint union of countably many connected open sets and $u_0$
is constant on every connected open subset of $U \cup W$, the set
\[\Theta:=\{u_0(x): \ \ x \in \overline{U \cup W} \}\]
is countable.

Now suppose $u_1$ is another minimizer. Then we have
\[\nabla u_0=0 \ \ \hbox{in} \ \ U \ \ \hbox{and}\ \ \frac{\partial u_0}{\partial \nu}=0 \ \ \hbox{on} \ \ \partial (V\cup W).\]

 Without loss of
generality we can assume $u_0\geq0$ in $\overline{\Omega}$. Then
\begin{eqnarray}\label{pro}
F(u_1)&=&\int_{\Omega  \setminus \overline{U\cup V\cup W}}
\sigma|\nabla u_0|.|\nabla u_1|dx \geq
\int_{\Omega  \setminus \overline{U\cup V\cup W}} \sigma|\nabla u_0.\nabla u_1|dx \nonumber \\
&\geq& \int_{\Omega  \setminus \overline{U\cup V\cup W}}\sigma
\nabla u_0.\nabla u_1=\int_{\partial \Omega}\sigma_0\frac{\partial
u_0}{\partial \nu}u_1 ds =\int_{\partial \Omega}\sigma_0\frac{\partial
u_0}{\partial \nu}f ds\\
&=&F(u_0),\nonumber
\end{eqnarray}
where $\nu$ is the outer normal to the boundary of $\Omega$. Since $u_0$ and $u_1$ both minimize the functional $F(u)$, equality
holds in (\ref{pro}). On the other hand the equality in
Cauchy inequality can only hold for parallel vectors, we have that
\begin{equation}
\nabla u_1(x)=\lambda(x)\nabla u_0(x), \ \ \ \ \hbox{a.e.} \ \ x \in
\Omega \setminus \overline{U\cup V\cup W},
\end{equation}
for some Lebesgue- measurable $\lambda(x)$. In particular,
\begin{equation}\label{normal.equ1}
\frac{\nabla u_0(x)}{|\nabla u_0(x)|}=\frac{\nabla u_1(x)}{|\nabla
u_1(x)|}
\end{equation} a.e. on
\[(\Omega \setminus \overline{U\cup V\cup W} )\cap \{x \in \Omega: |\nabla u_1|\neq 0\}.\]

Let $E_t=\{x \in \Omega \setminus \overline{U\cup V\cup W}:
u_0(x)>t\}$. Since $\Theta$ is countable, for a.e. $t>0$, $\partial
E_t \cap \overline{(U\cup W)}=\emptyset$ (otherwise $u_0$ must be a
constant). We claim that the sets $\partial E_t \cap (\Omega
\setminus \overline{V}) $ are smooth $C^1$ manifolds in $\Omega
\setminus \overline{V}$ for almost all $t>0$ with $\partial E_t \cap
\overline{U\cup W}=\emptyset$. To prove this note that since $u_0
\in C^{1}(\Omega \backslash \overline{U\cup V)}$, from equality
(\ref{normal.equ1}) we have that the measure theoretical normal
$\nu_{t}(x)=-\frac{\nabla u_0}{|\nabla u_0|}$ extends continuously
from $\partial^{*}E_t \cap (\Omega \backslash \overline{V})$ to the
topological boundary $\partial E_t \cap (\Omega \backslash
\overline{V})$, where $\partial^* E_t$ is the measure theoretical
boundary of $E_t$. By the regularity result of De Giorgi (see, e.g.
Theorem 4.11 in \cite{G84}), we conclude that $\partial E_t \cap
\Omega \backslash \overline{V}$ is a $C^1$-hypersurface for almost
all $t>0$.

The function $u_1$ is constant on every $C^1$ connected components
of $\partial E_t \cap (\Omega \backslash \overline{V})$. Indeed, let
$\gamma:(-\epsilon,+\epsilon) \rightarrow
\partial E_t \cap (\Omega \backslash \overline{V})$ be an arbitrary $C^1$ curve in $\partial E_t \cap (\Omega \backslash \overline{V})$. Then we
have

\[\frac{d}{dt}u_1(\gamma(s))=|\nabla u_1(\gamma(s))|\nu(\gamma (s)).\gamma'(s)=0,\]
because either $|\nabla u_1(\gamma(s))|=0$ or $\nu(\gamma
(s)).\gamma'(s)=0$ on $\partial E_t \cap (\Omega \backslash
\overline{V})$. So $u_1$ is constant along $\gamma$.

Let $t$ be one of the values for which $\partial E_t \cap (\Omega
\backslash \overline{V})$ is a hypersurface and $\partial E_t \cap
\overline{U \cup W}=\emptyset$ (which is the case for almost every
$t>0$). We show next that each connected component of $\partial E_t$
intersects the boundary $\partial \Omega$.

Arguing by contradiction, assume that $\Sigma_t$ is a connected
component of $\partial E_t$ such that $\Sigma_t \cap
\partial \Omega=\emptyset$. We consider two cases:\\

(I) $\Sigma_t \cap \partial V= \emptyset$,

(II) $\Sigma_t \cap
\partial V\neq \emptyset.$\\

Case I: Assume that $\Sigma_t \cap \partial V= \emptyset$. Then
$\partial \Omega \cup \Sigma_t$ is a compact manifold with two
connected components. By the Alexander duality theorem for $\partial \Omega \cup \Sigma_t$ (see, e.g., Theorem
27.10 in \cite{G81}) we have that $\R^n \setminus(\partial \Omega
\cup \Sigma_t)$ is partitioned into three open connected components:
$\R^n=(\R^n\setminus \overline{\Omega} \cup O_1\cup O_2)$. Since
$\Sigma_t \subset \Omega$ we have $O_1 \cup O_2=\Omega \setminus
\Sigma_t$ and then $\partial O_i \subset
\partial \Omega \cup \Sigma_t$ for $i=1,2$.

We claim that at least one of the $\partial O_1$ or $\partial O_2$
is in $\Sigma_t$. Assume not, i.e. for $i=1,2$, $\partial O_i \cap
\partial \Omega\neq \emptyset$. Since $\partial \Omega$ is connected
(by assumption) we have that $O_1 \cup O_2\cup \partial \Omega$ is
connected which implies that $O_1 \cup O_2\cup(\R^n\setminus
\Omega)$ is also connected. Again by applying the Alexander duality
theorem for $\Sigma_t \subset \R^n$, we have that $\R^n \setminus
\Sigma_t$ has exactly two open connected components, one of which is
unbounded: $\R^n \setminus \Sigma_t=O_{\infty}\cup O_0$. Since
$O_1\cup O_2 \cup (\R^n \backslash \Omega)$ is connected and
unbounded, we have that $O_1\cup O_2 \cup (\R^n \backslash \Omega)
\subset O_{\infty}$, which leaves $O_0 \subset \R^n \setminus (O_1
\cup O_2 \cup (\R^n \setminus \Omega)) \subset \Sigma_t$. This is
impossible since $O_0$ is open and $\Sigma_t$ is a hypersurface.
Therefore either $O_1$ or $O_2$ or both has the boundary in
$\Sigma_t$.

Assume $\partial O_1\subset \Sigma_t$. We claim that $u_0=t$ in
$O_1$. Indeed, since $O_1$ is an extension domain ($\partial
\Omega_1$ has a unit normal everywhere) the new map $\tilde{u_0}$
defined by

\begin{eqnarray*}
\tilde{u_0}:=\left\{ \begin{array}{ll}
u_0,\ \ x \in \Omega \setminus O_1,\\
t, \ \ \ \ x\in \overline{O_1},
\end{array} \right.
\end{eqnarray*}
is in $W^{1,1}(\Omega) \cap C(\overline{\Omega})$ and decreases the
functional, which contradicts the minimality of $u_0$. Therefore
$u_0=t$ in $O_1$, which makes $|\nabla u_0|=0$ in $O_1$. This is
contradiction since we have assumed $\partial E_t \cap \overline{U
\cup W}=\emptyset$.

Case (II): Assume $\Sigma_t \cap \partial V\neq \emptyset$ and let

\[V_{t}=\{x: \ \ x\in V_i \ \ \hbox{and}\ \ V_i\cap \Sigma_t\neq \emptyset\},\]
where $V_i$ are the connected components of $V$. Now define

\[\Sigma^*_t:=\partial V_t \cup \Sigma_t.\]
By our assumptions $\Sigma^*_t$ is a piecewise $C^1$-hyperfurface
and $\Sigma^*_t \cap \partial \Omega =\emptyset.$ Since $\partial
\Omega \cup \Sigma^*_t$ is a compact manifold with two connected
components, by the Alexander duality theorem and an argument similar to
that of case (I) we conclude that $\Omega \backslash \Sigma^*_t=O_1
\cup O_2$ and at least one of the $\partial O_1$ or $\partial O_2$
is in $\Sigma^*_t$. Assume $\partial O_1 \subset \Sigma^*_t$ and let
\[O=O_1 \cap (\Omega \backslash \overline{V}).\]
Then $O$ is a non-empty open subset of $\Omega \backslash
\overline{V}$. We claim that $u_0=t$ in $O$. Indeed the new map
defined by
\begin{eqnarray*}
\tilde{u_0}:=\left\{ \begin{array}{ll}
u_0,\ \ x \in \Omega \backslash (V \cup O),\\
t, \ \ \ \ x\in \overline{O},
\end{array} \right.
\end{eqnarray*}
can be extended to a function in $W^{1,1}(\Omega) \cap
C(\bar{\Omega})$ which decreases the functional and contradics the
minimality of $u_0$. Hence $u_0=t$ in $O$ which is a contradiction
because we have assumed $E_t \cap \overline{U\cup W}=\emptyset.$

In both cases the contradiction follows from the assumption that
$\Sigma_t \cap
\partial \Omega=\emptyset$. We conclude that each connected
component of $\partial E_t$ reaches the boundary $\partial
\Omega_t$. Since $u_0$ and $u_1$ coincide on the boundary $\partial
\Omega$, we have showed that $u_0|_{\partial E_t}=u_1|_{\partial
E_t}=t$ for almost every $t$. Therefore $u_0=u_1$ a.e. in
$\Omega\setminus \overline{U\cup W}$.


Now note that $u_0=u_1$ on the boundary of each connected component
of $U\cup W$. Since, $u_0$ and $u_1$ are constant on each connected
component of $U\cup W$, $u_0$ and $u_1$ should also agree on
$U\cup W$. Hence $u_0=u_1$ on $\Omega\setminus V$ and the proof is complete.  \hfill $\Box$ \\

%
{\bf Proof of Theorem \ref{p.unique}:}  To prove the theorem we
shall prove the stronger statement 2). It is enough to prove the
theorem for each connected component of $\Omega_{\alpha, \beta}$.
Hence without loss of generality we may assume that
$\Omega_{\alpha,\beta}$ is connected. By the definition of
$\Omega_{\alpha,\beta}$ we have
\begin{equation}\label{bound} u_1(\partial \Omega_{\alpha,\beta}\setminus \Gamma)\subset
\{\alpha, \beta\}.\end{equation} Let $J(x):=J_1(x)=J_2(x)$ for $x
\in \Omega_{\alpha,\beta}$. By our assumptions $|J|>0$ a.e. in
$\Omega \setminus_{\alpha,\beta} \overline{U^1\cup V^1\cup W^1}$.
Hence, $|\nabla u_1|>0$ a.e. on $\Omega_{\alpha,\beta} \setminus
\overline{U^1\cup V^1\cup W}$. Since $U^1\cup W^1$ is a disjoint
union of countably many connected open sets and $u_1$ is constant on
every connected open subset of $U^1 \cup W$, the set
\[\Theta:=\{u_1(x): \ \ x \in \overline{U^1\cup W^1} \}\]
is countable. Without loss of generality we can assume $u_1\geq0$ in
$\Omega_{\alpha,\beta}$.

Since $J_1=J_2$ in $\Omega_{\alpha,\beta}$, we have that
\begin{equation}
\nabla u_1(x)=\lambda(x)\nabla u_2(x), \ \ \ \ \hbox{a.e.} \ \ x \in
\Omega_{\alpha,\beta}\backslash \overline{U^1\cup V^1\cup W^1},
\end{equation}
for some nonnegative Lebesgue-measurable function $\lambda(x)$. In
particular, for a.e. $ x \in \Omega_{\alpha,\beta} \setminus
\overline{U^1 \cup V^1\cup W^1}$ we must have
\begin{equation}\label{normal.equ}
\frac{\nabla u_1(x)}{|\nabla u_1(x)|}=\frac{\nabla u_2(x)}{|\nabla
u_2(x)|}.
\end{equation}

Let $E_t=\{x \in \Omega_{\alpha,\beta} \setminus \overline{U^1 \cup
V^1\cup W^1}: u_1(x)>t\}$. Since $\Theta$ is countable, for a.e.
$t>0$, $\partial E_t \cap \overline{U^1\cup W^1}=\emptyset$
(otherwise $u_1$ must be a constant). With an argument similar to
that of Theorem \ref{unique}, one can show that the sets $\partial
E_t \cap (\Omega_{\alpha,\beta}\backslash \overline{V^1})$ are
smooth $C^1$ manifolds in $\Omega_{\alpha,\beta}$ for almost all
$t>0$ with $\partial E_t \cap \overline{U^1\cup W^1}=\emptyset$ and
the function $u_2$ is constant on each connected components of
$\partial E_t \cap (\Omega_{\alpha,\beta}\backslash
\overline{V^1})$.

Now let $t\neq \alpha, \beta$ to be one of the values for which
$\partial E_t \cap (\Omega_{\alpha,\beta}\backslash \overline{V^1})$
is a hypersurface and $\partial E_t \cap \overline{U^1\cup
W^1}=\emptyset$ (which is the case for almost every $t>0$). We next
show that each connected component of $\partial E_t$ intersects
$\Gamma$.

Arguing by contradiction, assume that $\Sigma_t \subset
\Omega_{\alpha,\beta}$ is a connected component of $\partial E_t$
such that $\Sigma_t \cap
\partial \Omega=\emptyset$. We consider two cases:\\

(I) $\Sigma_t \cap \partial V^1= \emptyset$,

(II) $\Sigma_t \cap
\partial V^1\neq \emptyset.$\\

Case I: Assume $\Sigma_t \cap \partial V^1= \emptyset$. Then
$\partial \Omega \cup \Sigma_t$ is a compact manifold with two
connected components. By the Alexander duality theorem we have that
$\R^n \setminus(\partial \Omega \cup \Sigma_t)$ is partitioned into
three open connected components: $\R^n=((\R^n\setminus
\overline{\Omega}) \cup O_1\cup O_2)$. Since $\Sigma_t \subset
\Omega$ we have $O_1 \cup O_2=\Omega \setminus \Sigma_t$ and then
$\partial O_i \subset
\partial \Omega \cup \Sigma_t$ for $i=1,2$. With an argument similar to the one provided for the proof of
Theorem \ref{unique}, we can show that at least one of the $\partial
O_1$ or $\partial O_2$ is in $\Sigma_t$. Assume $\partial O_1\subset
\Sigma_t$. Since $u_1$ satisfies the elliptic equation
\[\nabla.(\sigma_1\nabla u_1)=0, \ \ \hbox{in} \ \ O_1\]
and $u_1=t$ on $\partial O_1$, $u_1=t$ in $O_1$ and therefor $|J|=0$
on $O_1$. This is a contradiction since we have assumed $\partial
E_t \cap \overline{U^1\cup W^1}=\emptyset$.

Case (II): Assume $\Sigma_t \cap \partial V^1\neq \emptyset$ and let

\[V_{t}=\{x: \ \ x\in V_i^1 \ \ \hbox{and}\ \ V_i\cap \Sigma_t\neq \emptyset\},\]
where $V_i^1$ are the connected components of $V^1$. Now define

\[\Sigma^*_t:=\partial V_t \cup \Sigma_t.\]
By our assumptions $\Sigma^*_t$ is a piecewise $C^1$-hyperfurface
and $\Sigma^*_t \cap \partial \Omega =\emptyset.$ Since $\partial
\Omega \cup \Sigma^*_t$ is a compact manifold with two connected
components, by Alexander duality theorem and an argument similar to
that of Theorem \ref{unique} we conclude that $\Omega \backslash
\Sigma^*_t=O_1 \cup O_2$ and at least one of the $\partial O_1$ or
$\partial O_2$ is in $\Sigma^*_t$. Assume $\partial O_1 \subset
\Sigma^*_t$ and let
\[O=O_1 \cap (\Omega_{\alpha, \beta} \backslash \overline{V^1}).\]
Then $O$ is a non-empty open subset of $\Omega_{\alpha,\beta}
\backslash \overline{V^1}$. We claim that $u_0=t$ in $O$. Indeed the
new map defined by
\begin{eqnarray*}
\tilde{u_0}:=\left\{ \begin{array}{ll}
u_0,\ \ x \in \Omega \backslash (V^1 \cup O),\\
t, \ \ \ \ x\in \overline{O},
\end{array} \right.
\end{eqnarray*}
can be extended to a function in $W^{1,1}(\Omega) \cap
C(\bar{\Omega})$ that solves the equation (\ref{pde_inclusions}).
Since the equation (\ref{pde_inclusions}) has a unique solution
$u=\tilde{u}$. Thus $u_0=t$ in $O$ which is a contradiction since we
have assumed $\partial E_t \cap \overline{U^1\cup W^1}=\emptyset$.

In both cases the contradiction follows from the assumption
$\Sigma_t \cap
\partial \Omega_{\alpha,\beta}=\emptyset$. Since $t\neq \alpha, \beta$ and
\[u_1(\partial \Omega_{\alpha,\beta}\setminus \Gamma)\subset
\{\alpha, \beta\},\] $E_t$ intersects $\Gamma$ for almost every
$t\geq 0$.

 Since $u_0$
and $u_1$ coincide on  $\Gamma$, we have showed that $u_1|_{\partial
E_t}=u_2|_{\partial E_t}=t$ for almost every $t$. Therefore
$u_0=u_1$ a.e. in $\Omega_{\alpha, \beta}\setminus \overline{U^1\cup
W^1}$. Now note that $u_1=u_2$ on the boundary of each connected
component of the set $U^1\cup W$. Since, $u_1$ and $u_2$ are
constant on each connected component of $U^1\cup W^1$, $u_1$ and
$u_2$ should also agree on $U^1\cup W$. Hence $u_1=u_2$ on
$\Omega_{\alpha,\beta}\backslash \overline{V^1}$. The proof is
complete. \hfill $\Box$

\section{Equipotential surfaces are area minimizing in the conformal metric}
In this section we present the proof of Theorem \ref{area.min}.  We
prove that the equipotential sets are global minimizers of
$E(\Sigma)$. This is a consequence of minimizing property of the
voltage potential for the functional $F(u)$. First we recall the
co-area formula.

\begin{theorem}\label{co-area} (Co-area formula). Let $u \in
Lip(\Omega)$ and $a$ be integrable in $\Omega \subset \R^n$. Then,
for a.e. $t \in \R$, $H^{n-1}(u^{-1}(t)\cap \Omega)<\infty$ and
\begin{equation}
\int_{\Omega} a|\nabla
u(x)|dx=\int_{-\infty}^{\infty}\int_{u^{-1}(t)}a dH^{n-1}(x)dt,
\end{equation}
where $H^{n-1}$ is the $(n-1)-$dimensional Hausdorff measure.\\
\end{theorem}

\begin{proposition}\label{pro1} Let $a\geq
0$ be integrable in $\Omega$, $U$ be an open subset of $\Omega$, and
\small{
\begin{equation*}
u\in argmin \left\{\int_{\Omega}a |\nabla v|dx: \ \ v \in
Lip(\Omega), \ \ \hbox{and }\ \ v|_{\Omega}=f \right\}.
\end{equation*}
} For $\lambda \in \R$ arbitrary fixed, let
$u_{+}=\max\{u-\lambda,0\}$ and $u_{-}=\max\{u,\lambda\}$ be defined
in $\Omega$, and $f_{+}=\max\{f-\lambda,0\}$, respectively
$f_{-}=\min \{f,\lambda\}$, be defined on the boundary $\partial
\Omega$.  Then \small{
\begin{equation*}
u_{+}\in argmin \left\{\int_{\Omega}a |\nabla v|dx: \ \ v \in
Lip(\Omega), \ \ \hbox{and }\ \ v|_{\Omega}=f_{+} \right\},
\end{equation*}
and
\begin{equation*}
u_{-}\in argmin \left\{\int_{\Omega}a |\nabla v|dx: \ \ v \in
Lip(\Omega),\ \ \hbox{and }\ \ v|_{\Omega}=f_{-} \right\}.
\end{equation*}
}
\end{proposition}
{\bf Proof:} The proof is similar to the proof of Proposition 2.2
\cite{NTT10} and we omit it. \hfill $\Box$\\

\begin{corollary}\label{cor1}
Let $a\geq 0$ be integrable in $\Omega$, $U$ be an open subset of
$\Omega$, and

\small{
\begin{equation*}
u\in argmin \left\{\int_{\Omega}a |\nabla v|dx: \ \ v \in
Lip(\Omega), \ \ \hbox{and }\ \ v|_{\Omega}=f \right\}.
\end{equation*}}
For every $\lambda \in \R$ and $\epsilon>0$ define
\begin{equation}\label{u.ep}
u_{\lambda,\epsilon}:=\frac{1}{\epsilon}\min\{\epsilon,\max\{u-\lambda,0\},\}
\end{equation}
and let $f_{\lambda,\epsilon}$ be its trace on the boundary
$\partial \Omega$. Then $u_{\lambda,\epsilon} \in Lip(\Omega)$ and
\begin{equation*}
u_{\lambda,\epsilon}\in argmin \left\{\int_{\Omega}a |\nabla v|dx: \
\ v \in Lip(\Omega), \ \ \hbox{and}\ \
v|_{\Omega}=f_{\lambda,\epsilon} \right\}.
\end{equation*}
\end{corollary}
{\bf Proof:} The proof follows directly from Proposition \ref{pro1}
applied twice. \hfill $\Box$\\
\begin{lemma}\label{lem}
Let $a, u\in  Lip(\Omega)$ such that
\[\{x: \ \ |\nabla u(x)|=0\}=Z\cup L,\]
where $Z$ is open and $L$ has Lebesgue measure zero,
$a(\overline{Z})=\{0\}$, and
\begin{equation}
a\frac{\nabla u}{|\nabla u|}\in W^{1,1}(\Omega \backslash
\overline{Z}).
\end{equation}
Then for almost every $\lambda \in \R$,
\begin{equation}\label{co-co2}
\lim_{\epsilon \rightarrow 0}\int_{\Omega}a|\nabla
u_{\lambda,\epsilon}| dx =\int_{u^{-1}(\lambda)}a dH^{n-1}(x),
\end{equation}
where $u_{\lambda,\epsilon}$ is defined by $(\ref{u.ep})$.
\end{lemma}
{\bf Proof.} The proof is similar to the proof of Lemma 2.4 in
\cite{NTT10}. From Theorem \ref{co-area}, we have
$H^{n-1}(u^{-1}(\lambda)\cap \Omega)<\infty$, a.e. $\lambda \in \R$.
In particular
\begin{equation}\label{n.m }
H^{n}(u^{-1}(\lambda)\cap \Omega)=0.
\end{equation}
Since $H^{n-1}(\partial \Omega)<\infty$, from the disjoint partition
$\partial \Omega=\bigcup_{\lambda \in \R}(u^{-1}(\lambda)\cap
\partial \Omega)$, we have
\[H^{n-1}(u^{-1}(\lambda)\cap \partial \Omega)>0,\]
for at most countable many $\lambda$. In particular, for almost
every $\lambda \in \R$
\begin{equation}\label{(n-1).m}
H^{n-1}(u^{-1}(\lambda)\cap
\partial \Omega)=0.
\end{equation}
Let $\lambda \in Range (u)$ be such that both (\ref{n.m }) and
(\ref{(n-1).m}) hold, and $\epsilon>0$. Recall

\begin{eqnarray*}
u_{\lambda,\epsilon}=\left\{ \begin{array}{ll} 0 &\hbox{if}\ \
u(x)<\lambda,\\
(u(x)-\lambda)/\epsilon & \hbox{if} \ \ \lambda\leq u(x) \leq
\lambda +\epsilon,\\
0  & \hbox{if} \ \ u(x)>\lambda+\epsilon.
\end{array} \right.
\end{eqnarray*}
From the co-area formula we have
\begin{eqnarray}\label{co-co1}
\int_{\Omega}a|\nabla u_{\lambda,\epsilon}|dx &=&\int_{-\infty}^{+
\infty}\int_{(u_{\lambda,\epsilon})^{-1}(t)}a d H^{n-1}(x)dt\nonumber \\
&=&\int_{0}^{1}\int_{\{x:\ \ u(x)=\lambda+t\epsilon\}}a d H^{n-1}(x)
\end{eqnarray}
To complete the proof it is enough to prove that
\begin{equation}\label{Co-co}
\lim_{\epsilon \rightarrow 0}\int_{\{x: u(x)=\lambda+\epsilon\}}a
dH^{n-1}(x)=\int_{\{x: u(x)=\lambda\}}a dH^{n-1}(x),
\end{equation}
holds uniformly for almost every $t\in [0,1]$. The domain
\[\Omega_{t,\epsilon}:=\{x\in \Omega: \lambda<u(x)<\lambda+t\epsilon\}\]
is Lipschitz. Since $a \in Lip(\Omega)$, it extends continuously to
the boundary. The $a\nabla u/|\nabla u| \in W^{1,1}(\Omega
\backslash \overline{Z})$ also extends to the boundary $\partial
(\Omega \backslash Z)$ as a bounded function. Now notice that
$u(\overline{Z})$ is at most countable. Therefore, for a. e.
$\lambda \in Range(u)$ and a.e. $t\in[0,1]$ the outer unit normal
$\nu$ to the boundary $\partial \Omega_{t,\epsilon}$ exists. Then
Green's formula in $\Omega_{t,\epsilon}$ yields \small{
\begin{eqnarray*}
& & |\int_{u^{-1}(\lambda+t\epsilon)}a d
H^{n-1}-\int_{u^{-1}(\lambda)} a d H^{n-1}|\\&=&
|\int_{u^{-1}(\lambda+t\epsilon)\cap \Omega \backslash
{\overline{Z}}}a d H^{n-1}-\int_{u^{-1}(\lambda) \cap \Omega
\backslash {\overline{Z}}} a d H^{n-1} |\\
& \leq &|\int_{\{x \in \partial \Omega: \ \
\lambda<u(x)<\lambda+\epsilon\}}a\frac{\nabla u}{|\nabla u|}\nu d
H^{n-1}(x)|\\
&+&|\int_{\{x \in  \Omega: \ \
\lambda<u(x)<\lambda+\epsilon\}}\nabla. a\frac{\nabla u}{|\nabla
u|}dx|.
\end{eqnarray*}

} Using (\ref{n.m }) we have
\begin{eqnarray*}
\lim_{\epsilon \rightarrow 0} H^n(\{x\in \Omega : \lambda
<u(x)<\lambda+\epsilon\})&=&H^n( \bigcap_{\epsilon>0}\{x\in \Omega :
\lambda <u(x)<\lambda+\epsilon\})\\
&\leq& H^n( \bigcap_{\epsilon>0}\{x\in \Omega : \lambda \leq
u(x)<\lambda+\epsilon\})\\
&=&H^n(u^{-1}(\lambda)\cap \Omega)=0.
\end{eqnarray*}
Similarly by (\ref{(n-1).m}) we obtain
\[\lim_{\epsilon \rightarrow 0}H^{n-1}(\{x\in \partial \Omega : \lambda
< u(x)<\lambda+\epsilon\})=0.\] This proves (\ref{Co-co}). By taking
the limit $\epsilon \rightarrow 0$ in (\ref{co-co1}) and using
(\ref{Co-co}) we obtain (\ref{co-co2}). \hfill $\Box$\\

{\bf Proof of Theorem \ref{area.min}:} For $\lambda \not \in
Range(u)$, the left hand side of (\ref{ine}) is zero and and the
inequality trivially holds. Since $u$ obeys the maximum principle
and $u=v$ on $\partial \Omega$, $Range(u)\subset Range(v)$.

Now let $\lambda \in Range(u) \setminus \left(u(\overline{Z})\cup
v(\overline{Z_{v}}) \right)$ and recall that $u(\overline{Z})$ and
$u(\overline{Z_v})$ are both countable. Since $|\nabla u|\neq 0$
a.e. in $\Omega \setminus \overline{Z}$ and $|\nabla v|\neq 0$ a.e.
in $\Omega \setminus \overline{S_v}$, for almost every $\lambda \in
Range(u)$ the corresponding $\lambda-$level set is a $C^1$-smooth
oriented surface. In particular the $H^{n-1}-$measure coincides with
the induced Lebesgue measure on the respective surface. Moreover,
$u$ and $v$ satisfy $(\ref{n.m })$ and $(\ref{(n-1).m})$ for a.e.
$\lambda \in \R$.

For $\epsilon>0$ arbitrary fixed, let $u_{\lambda,\epsilon}$ be
defined by (\ref{u.ep}) and define similarly
\[v_{\lambda,\epsilon}:=\min\{\epsilon, \max\{v-\lambda,0\}/\epsilon\}.\]
Since $u=v$ on the boundary $\partial \Omega$, we also have
$u_{\lambda,\epsilon}=v_{\lambda,\epsilon}$ on $\partial \Omega$.
From Corollary \ref{cor1} we have
\begin{equation}
\int_{\Omega}a |\nabla u_{\lambda,\epsilon}|dx \leq \int_{\Omega}a
|\nabla v_{\lambda,\epsilon}|dx.
\end{equation}
Letting $\epsilon \rightarrow 0$ and applying Lemma \ref{lem} we
obtain (\ref{ine}). \hfill $\Box$\\

\section{Appendix: Perfectly conductive and insulating inclusions}
The results in this appendix formalize the definition of perfectly
conducting  as infinity limit of
conductivity. They are slight generalization of the ones in
\cite{baoLiYin} to include both perfectly conductive and insulating
inclusions.

Let $U=\cup_{j=1}^\infty U_j$ be an open subset of $\Omega$ with
$\overline{U}\subset\Omega$ to model the union of the connected
components  $U_j$ ($j=1,2,...$) of \emph{perfectly conductive
inclusions}, and $V$ be an open subset of $\Omega$ with
$\overline{V}\subset\Omega$ to model the union of all connected
\emph{insulating inclusions}. Let
 $\chi_U$ and $\chi_V$ be their corresponding characteristic
function. We assume that $\overline{U}\cap\overline{V}=\emptyset$,
 $\Omega\setminus \overline{U\cup V}$ is connected, and that the
boundaries $\partial U$, $\partial V$ are piecewise $C^{1,\alpha}$.
Let $\sigma_1\in L^\infty(U)$, and $\sigma\in
L^\infty(\Omega\setminus \overline{U\cup V})$ be such that
\begin{equation}\label{ellipticity}
0<\lambda\leq\sigma_1,\sigma\leq \Lambda<\infty,
\end{equation}for some positive constants $\lambda$ and $\Lambda$.

For each $0<k<1$ consider the conductivity problem
\begin{equation}
\nabla\cdot(\chi_U(\frac{1}{k} \sigma_1-\sigma)+\sigma)\nabla u=0,
\qquad \frac{\partial u}{\partial \nu}=0 \ \hbox{on} \ \
\partial V, \ \ \hbox{and} \qquad u|_{\partial\Omega}=f.\label{condEQ}
\end{equation} 
The condition on $\partial V$ ensures that $V$ is insulating. 
It is well known that the problem (\ref{condEQ}) has a unique
solution $u_k\in H^1(\Omega)$ which also solves:
\begin{equation}
\left\{ \begin{array}{ll}
\nabla\cdot \sigma \nabla u_k=0,&\mbox{in}\,\Omega\setminus\overline{U\cup V},\\
\nabla\cdot \sigma_1 \nabla u_k=0,&\mbox{in}\,U,\\
u_k|_+=u_k|_-,&\mbox{on}\,\partial U,\\
\left. \frac{1}{k} \sigma_1\frac{\partial u_k}{\partial \nu}\right|_-=\left.\sigma\frac{\partial u_k}{\partial \nu}\right|_+,&\mbox{on}\,\partial U,\\
\left. \frac{\partial u_k}{\partial \nu}\right|_{+}=0,&\mbox{on}\,\partial V,\\
u_k|_{\partial\Omega}=f.
\end{array}\label{pde_k} \right.
\end{equation} Moreover, the energy functional
\begin{equation}
I_k[v]=\frac{1}{2k}\int_U\sigma_1|\nabla
v|^2dx+\frac{1}{2}\int_{\Omega\setminus \overline{U\cup
V}}\sigma|\nabla v|^2 dx
\end{equation}
has a unique minimizer over the maps in $H^1(\Omega)$ with trace $f$
at $\partial\Omega$ which is the unique solution $u_k$ of
(\ref{pde_k}).

We shall show below why the limiting solution (with $k\to 0$) solves

\begin{equation}\left\{ \begin{array}{ll}
\nabla\cdot \sigma \nabla u_0=0,&\mbox{in}\,\Omega\setminus\overline{U\cup V},\\
\nabla u_0=0, &\mbox{in} \ \ U,\\
u_0|_+=u_0|_-,&\mbox{on}\ \ \partial U,\\
\int_{\partial U_j}\sigma\frac{\partial u_0}{\partial \nu}|_{+}ds=0,&j=1,2,...,\\
\frac{\partial u_0}{\partial \nu}|_{+}=0,&\mbox{on}\;\partial V,\\
u_0|_{\partial\Omega}=f,\\
\end{array}\label{pde_0} \right.
\end{equation}
By elliptic regularity $u_0 \in C^{1,\alpha}(\Omega \backslash U\cup
V)$ and for any $C^{1,\alpha}$ boundary portion $T$ of $\partial
(U\cup V)$, $u_0 \in C^{1,\alpha}((\Omega \backslash (U\cup V))\cup
T)$.

\begin{proposition}
The problem (\ref{pde_0}) has a unique solution in $H^1(\Omega)$
which is the unique minimizer of the functional
\begin{equation}\label{I0}
I_0[v]=\frac{1}{2}\int_{\Omega\setminus \overline{U\cup
V}}\sigma|\nabla v|^2 dx,
\end{equation}over the set $A_0:=\{u\in H^1(\Omega\setminus\overline{V});\,
u|_{\partial\Omega} =f,\, \nabla u=0\, \mbox{in}\, U\}$.
\end{proposition}
{\bf Proof:} Note that $A_0$ is weakly closed in
$H^1(\Omega\setminus \overline{V})$. The functional $I_0$ is lower
semicontinuous, strictly convex, and, thus, has a unique minimizer
$u_0^*$ in $A_0$.

First we show that $u_0^*$ is a solution of (\ref{pde_0}). Since
$u_0^*$ minimizes (\ref{I0}), we have
\begin{equation}\label{weaksolution}
0=\int_{\Omega\setminus \overline{V\cup U}}\sigma \nabla
u^*_0\cdot\nabla\phi dx,
\end{equation}
for all $\phi\in H^1(\Omega\setminus \bar{V})$, with
$\phi|_{\partial\Omega}=0$, and $\nabla \phi =0 $ in $U$. In
particular, if $\phi\in H^1_0(\Omega\setminus \bar{V})$, we get
$\int_{\Omega\setminus \overline{U\cup V}}(\nabla \cdot\sigma\nabla
u_0^*)\phi dx=0$ and thus $u_0^*$ solves the conductivity equation
in (\ref{pde_0}). If we choose $\phi\in H^1(\Omega\setminus
\bar{V})$, with $\phi|_{\partial\Omega}=0$, and $\phi \equiv 0 $ in
$U$, from Green's formula applied to (\ref{weaksolution}), we get
$\int_{\partial V} \sigma\left.\frac{\partial
u_0^*}{\partial\nu}\right|_+\phi=0, \,\forall \,\phi|_{\partial
V}\in H^{1/2}(\partial V),$ or, equivalently,
$\sigma\left.\frac{\partial u_0^*}{\partial \nu}\right|_{\partial
V}=0$. If we choose $\phi_j\in H^1_0(\Omega\setminus \bar{V})$ with
$\phi_j\equiv 1$ in the connected component $U_j$ of $U$ and
$\phi_j\equiv 0$ in $U\setminus U_j$, from Green's formula applied
to (\ref{weaksolution}) we obtain $\int_{\partial
U_j}\sigma\frac{\partial u_0^*}{\partial\nu}=0$.

Next we show that the equation (\ref{pde_0}) has a unique solution
and, consequently, $u_0^*=u_0|_{\Omega\setminus\overline{V}}$.
Assume that $u^1$ and $u^2$ are two solutions and let $u=u_2-u_1$,
then $u|_{\partial\Omega}=0$ and

\begin{eqnarray}0=&-\int_{\Omega\setminus\overline{U\cup V}}(\nabla
\cdot\sigma\nabla u)udx=-\int_{\partial\Omega}\sigma\frac{\partial
u}{\partial \nu}uds+\int_{\partial V} \sigma\left.\frac{\partial
u}{\partial \nu}\right|_+uds\\&+\int_{\partial
U}\sigma\left.\frac{\partial u}{\partial \nu}\right|_+uds
+\int_{\Omega\setminus\overline{U\cup V}}\sigma|\nabla
u|^2dx=\int_{\Omega\setminus\overline{U\cup V}}\sigma|\nabla u|^2dx.
\end{eqnarray}
Since $\sigma\geq\lambda>0$, we get $|\nabla u|=0$ in
$\Omega\setminus\overline{V}$. Since $\Omega\setminus\overline{V}$
is connected and $u=0$ at the boundary, we conclude uniqueness of
the solution of the equations (\ref{pde_0}). \hfill $\Box$
\begin{theorem}
Let $u_k$ and $u_0$ be the unique solution of (\ref{pde_k})
respectively (\ref{pde_0}) in $H^1(\Omega)$. Then
$u_k\rightharpoonup u_0$ and, consequently, $I_k[u_k]\to I_0[u_0]$
as $k\to 0^+$.
\end{theorem}

{\bf Proof:} We show first that $\{u_k\}$ is bounded in
$H^1(\Omega)$ uniformly in $k\in (0,1)$. 
%
Since $1/k>1$, we have
\begin{eqnarray*}
\frac{\lambda}{2}\|\nabla
u_k\|^2_{L^2(\Omega\setminus\overline{V})}\leq\frac{1}{2}\int_{\Omega\setminus\overline{U\cup
V}}\sigma|\nabla u_k|^2dx+\frac{1}{2k}\int_U\sigma_2|\nabla
u_k|^2dx\\ \leq I_k[u_k]\leq I_k[u_0]\leq\frac{\Lambda}{2}\|\nabla
u_0\|^2_{L^2(\Omega\setminus\overline{V})},
\end{eqnarray*}or
\begin{equation}\label{ineqb}
\|\nabla u_k\|^2_{L^2(\Omega\setminus
\overline{V})}\leq\frac{\Lambda}{\lambda}\|\nabla
u_0\|^2_{L^2(\Omega\setminus\overline{V})}.
\end{equation}
From (\ref{ineqb}) and the fact that $u_k|_{\partial\Omega}=f$, we
see that $\{u_k\}$ is uniformly bounded in $H^1(\Omega \backslash
V)$ and hence weakly compact. Therefore, on a subsequence
$u_k\rightharpoonup u_0^*$ in $H^1(\Omega \backslash V)$, for some
$u_0^*$ with trace $f$ at ${\partial\Omega}$.

We will show next that $u_0^*$ satisfies the equations
(\ref{pde_0}), and therefore $u_0^*=u_0$ on $\Omega$. By the
uniqueness of solutions of (\ref{pde_0}) we also conclude that the
whole sequence converges to $u_0$.

Since $u_k\rightharpoonup u_0^*$ we have that
$0=\int_{\Omega\setminus \overline{U\cup V}}\sigma\nabla
u_k\cdot\nabla \phi dx\to \int_{\Omega\setminus \overline{U\cup
V}}\sigma_2\nabla u_0^*\cdot\nabla\phi dx$, for all $\phi\in
C^\infty_0({\Omega\setminus \overline{U\cup V}})$. Therefore $\nabla
\cdot\sigma\nabla u_0^*=0$ in ${\Omega\setminus \overline{U\cup
V}}$. Also because $u_k$ is a  minimizers of $I[u_k]$ we must have
$\nabla u_0^*=0$ in $U$. To check the boundary conditions, note
that, for all $\phi\in C^\infty_0(\Omega)$ with $\phi\equiv 0$ in
$U$, we have $\int_{\partial V}\sigma\left.\frac{\partial
u_k}{\partial\nu}\right|_+\phi ds=0$. Using the fact that $\phi$
were arbitrary, by taking the weak limit in $k\to 0$, we get
$\left.\frac{\partial u_0^*}{\partial\nu}\right|_+=0$ on $\partial
V$.  A similar argument applied to $\phi\in C^\infty_0(\Omega)$ with
$\phi\equiv 0$ in $V$, $\phi\equiv 1 $ in $U_j$, and $\phi\equiv 0$
in $U \backslash U_j$, also shows that $\int_{\partial
{U_j}}\sigma\left.\frac{\partial u_0^*}{\partial\nu}\right|_+\phi
ds=0$. Hence $u^*_0$ is the unique solution of the equation
(\ref{pde_0}) on $\Omega \backslash \overline{V}$. Thus $u_k$
converges weakly to the solution $u_0$ of (\ref{pde_0}) in
$\Omega\backslash \overline{V}$. \hfill $\Box$\\

 \end{document}